\documentclass[12pt]{article}
\usepackage{amssymb,latexsym,amsmath}
\usepackage{amsthm}
\usepackage{geometry}
\usepackage{graphicx,color}
\usepackage{lineno}
\theoremstyle{plain}

\theoremstyle{definition}

\numberwithin {equation}{section}
\usepackage{longtable}

\usepackage{float}
\restylefloat{table}

\begin{document}

\title
{A dynamic game on Green Supply Chain Management}
\author{Mehrnoosh Khademi $^a$\thanks{mehrnushkhademi@yahoo.com}
\and Massimiliano Ferrara $^{b,c,d}$\thanks{massimiliano.ferrara@unirc.it} \and Bruno Pansera$^b$\thanks{Presenter:
bruno.pansera@unirc.it} \and Mehdi Salimi
$^{d,e}$\thanks{mehdi.salimi@tu-dresden.de \&
mehdi.salimi@medalics.org}}
\date{}
\maketitle
\begin{center}
$^{a}$Department of Industrial Engineering, Mazandaran University
of Science and Technology, Babol, Iran\\ $^{b}$Department of Law
and Economics - University Mediterranea of Reggio Calabria,
Italy\\ $^{c}$CRIOS - Center for Research in Innovation,
Organization and Strategy - Department of Management and
Technology - Bocconi University, Italy\\ $^{d}$Medalics, Research
Centre of the University for Foreigner "Dante Alighieri" in Reggio
Calabria, Italy\\ $^{e}$Center for Dynamics, Department of
Mathematics, Technische Universit{\"a}t Dresden, Germany
\end{center}
\maketitle


\begin{abstract}\noindent

In this paper, we establish a dynamic game to allocate CSR
(Corporate Social Responsibility) to the members of a supply
chain. We propose a model of three-tier supply chain in
decentralized state that is including supplier, manufacturer and
retailer. For analyzing supply chain performance in decentralized
state and the relationships between the members of supply chain,
we use Stackelberg game and we consider in this paper a
hierarchical equilibrium solution for a two-level game. Specially,
we formulate a model that crosses through multi-periods by a
dynamic discreet Stackelberg game. We try to obtain an equilibrium
point at where both the profits of members and the level of CSR
taken by supply chains are maximized.


 \end{abstract}

\section{Introduction}

In recent years, a growing number of large, medium, and even
small- sized companies have increasingly focused on CSR. They have
realized the need to develop strategies that extend their
traditional corporate governance processes beyond firm boundaries
to their supply chain partners \cite{Kytle}. This is chiefly
because, along with increasing consumer information about the
conditions of manufacture, they criticize supply chains for
several social responsibilities. Moreover, the firms in supply
chain have been pressured by a regulations and policies related to
CSR from governments and organizations.
 The members of supply
chain make their decisions based on maximizing of their individual
net benefits. Also, when they have to take a level of CSR; this
situation leads to an equilibrium status. Game theory is one of
the most effective tools to deal with such a kind of management
problems.

A growing number of research papers, use game theoretical
applications in supply chain management. Cachon et al.
\cite{Cachon} discuss Nash equilibrium in a noncooperative cases
in a supply chain where there are one supplier and multiple
retailers. Hennet et al. \cite{Hennet} presented a paper to
evaluate the efficiency
 of different types of contracts between the industrial partners of a supply chain.
They applied game theory method for decisional purposes.

Tian et al. \cite{Tian} presented a system dynamics model based on
evolutionary game theory for green supply chain management as
well.

In this paper, we formulate a model for decentralized supply chain
network in CSR conditions in a long term with one leader and two
followers. The Stackelberg game model is recommended and applied
here to find an equilibrium point in which we maximize the profit
of members of supply chain and the level of CSR taken by the
supply chain. In this research, the supplier as a leader, can know
the optimal reaction of his followers, and regards such processes
to maximize his own profit. The manufacturer and the retailer as
followers, try to maximize their profits by considering all
conditions. We propose a Hamiltonian matrix to solve the optimal
control problem to obtain the equilibrium in this game. The paper
is organized as follows: Section 2 is devoted to General model.
Objective functions, constraints and solving the game are
illustrated in Section 3. A Conclusion is provided in Section 4.

\section{The general model}
We consider a Stackelberg differential game involving two players
playing the game over a fixed finite horizon model suggested by He
et al. \cite{He}. We consider a dynamic game that goes through
multi-periods as a repeated game with complete information. This
model is a three-tier, multi-period, decentralized supply chain
network. We assume that only one supplier, one manufacturer and
one retailer are involved in playing the Stackelberg game as well
as allocated social responsibility. A long term Stackelberg game
is played between the members of the decentralized supply chain
through two levels in which all members take CSR into
consideration. We formulate the model by selecting the supplier as
the leader and both of the manufacturer and retailer as the
follower in Stackelberg game. This model can be solved by
considering two levels of Stackelberg game. In first level the
manufacturer as the leader and retailer as the follower are
considered. In this level, we find equilibrium point and in the
second level, we consider the supplier as the leader and the
manufacturer as the follower. In fact, we put response functions
of followers in the objective function of the leader and we find
the final equilibrium point and all of the players make decisions.

As any problem formulated as a dynamic game, this model has a
state variable and control variables. We define the state variable
as the level of social responsibility taken by companies, and the
control variables are the capital amounts invested in taking
social responsibility. Specifically, all of the social
responsibility taken by firm $j$ at period $t$ can be expressed as
investment $I_t^j$. The level of current supply chain investment
in supply responsibility is $x_t$; therefore the accumulation of
level of social responsibility taken by the firms is given by
$x_{t+1}=\alpha x_t+\beta_1I_t^S+\beta_2I_t^{M}+\beta_3I_t^{R}$.

Here, ${\beta_1}$ is the rate of converting the supplier's capital
investment in CSR to the amount of CSR taken by the supply chain;
${\beta_2}$ is the rate of converting the manufacturer's capital
investment in CSR to the amount of CSR taken by the supply chain
and ${\beta_3}$ is the rate of converting the retailer's capital
investment in CSR to the amount of CSR taken by the supply chain
\cite{Shi}.

For the purpose of the paper, we more specifically assume:
function $B_t(x_t)=\delta x_t$, that represents a social benefit
to the firms, where the coefficient $\delta$ is supposed
 to be strictly positive \cite{Batabyal}.

The following functions $T_t^S=\tau
I_t^S\big[1+\theta(I_t^S+I_t^M+I_t^R)\big]$, $T_t^M=\tau
I_t^M\big[1+\theta(I_t^S+I_t^M+I_t^R)\big]$ measure the amount of
money given by the to the supplier and the manufacturer
\cite{Feibel}. For retailer $T_t^R=\tau
I_t^R\big[1+\theta(I_t^S+I_t^M+I_t^R)\big]$ both $\tau$ and
$\theta$ are tax return policy parameters. Specifically, $\tau$ is
the rate of individual post tax return on investment (ROI), and
$\theta$ is rate of supply chain's post tax return on investment
(ROI).

All kinds of social responsibilities are assumed to be expressed
as investment {$I_t$}.

The market inverse demand is $P^M(q_t)=a-bq_t$ \cite{Mankiw}.

\section{Objective function and constraints}
Let the time interval be [1,T]. The objective function of the
supplier is

\begin{equation*}
\begin{split}
J^S&=\arg \max
\sum_{t=1}^TP_t^Sq_t-cq_t+B_t^S(x_t)+T_t^S(I_t^S,I_t)-I_t^S+dI_t^{M}\\
&=\arg \max \sum_{t=1}^T vq_t-cq_t+\delta x_t^2+\tau
I_t^S[1+\theta(I_t^S+I_t^{M}+I_t^{R})]-I_t^S+dI_t^{M},
\end{split}
\end{equation*}
subject to $x_{t+1}=\alpha
x_t+\beta_1I_t^S+\beta_2I_t^{M}+\beta_3I_t^{R}$, where the
coefficients $\beta$ are positive and with $\beta<1$.

$P_t^S$ is the price of the supplier's raw material. Let
${P_t^S}=v$. $B_t^{S}(x_t)$ is the social benefit of the supplier,
$\delta$ is the parameter of the supplier's social benefit and
$T_t^{S}(I_t^{S}, I_t)$ is the tax return of the supplier.  $d$ is
the percentage of investment of the supplier payoff. Similarly,
the objective function of the manufacturer is

\begin{equation*}
\begin{split}
J^M&=\arg \max
\sum_{t=1}^TP_t^M(q_t)q_t-P_t^Sq_t+B^M(x_t)+T^M(I_t^M, I_t)-I_t^M+\widehat{d}I_t^R\\
&=\arg \max
\sum_{t=1}^T(a-bq_t)q_t-vq_t+\widehat{\delta}x_t^2+\tau
I_t^M(1+\theta (I_t^S+I_t^M+I_t^R))-I_t^M+\widehat{d}I_t^R,
\end{split}
\end{equation*}
where $P_t^M(q_t)$ is the retail price of the product of the
manufacturer. $B_t^{M}(x_t)$ is the social benefit of the
manufacturer, $\widehat{\delta}$ is the parameter of the
manufacturer's social benefit. $T_t^{M}(I_t^{M}, I_t)$ is the tax
return of the manufacturer. $\widehat{d}$ is the percentage of
investment of the manufacturer payoff.\\
The objective function of the retailer is
\begin{equation*}
\begin{split}
J^R&=\arg \max
\sum_{t=1}^TP_t^Rq_t-P_t^M(q_t)q_t+B^R(x_t)+T^R(I_t^R, I_t)-I_t^R\\
&=\arg \max
\sum_{t=1}^Tzq_t-(a-bq_t)q_t+\widehat{\widehat{\delta}}x_t^2+\tau
I_t^R(1+\theta (I_t^S+I_t^M+I_t^R))-I_t^R,
\end{split}
\end{equation*}
where $P_t^{R}$ is the price of the product the retailer sells to
the consumer. Let $P_t^{R}=Z$. $B_t^{R}(x_t)$ is the social
benefit of the retailer, $\widehat{\widehat{\delta}}$ is the
parameter of the retailer's social benefit. $T_t^{R}(I_t^{R},
I_t)$ is the tax return of the retailer.

\subsection{ Mathematical model}

We solve the mathematical model with two levels, in the first
level the manufacturer's optimal function is calculated by
reaction function of retailer and the second level, the game is
between the supplier as the leader and the manufacturer as the
follower. In fact, the reaction functions of two followers
(retailer and manufacturer) are placed on the objective function
of the leader (supplier), we can find final equilibrium point.

In the level one, we establish a Stackelberg game between
manufacturer as the leader and retailer as the follower. In this
level, to calculate the equilibrium first we calculate the best
reaction function of retailer, then we determine the
manufacturer's optimal decisions based on the retailer' best
reactions. Since we consider this dynamic game as an optimal
control problem, the Hamiltonian function is a practical way to
find the equilibrium of the game \cite{Sethi}. The manufacturer's
optimal decisions based on the retailer' best reactions  is
determined. To obtain the Stackelberg strategy of the
manufacturer, we maximize the objective function of the
manufacturer by its Hamiltonian function.

We define the firm's Hamiltonian as below. For fixed $I_t^M$ the
Hamiltonian function of the retailer is defined by
\begin{equation*}
H_t^{R}=J_t^{R}+P_{t+1}^{R}( x_{t+1}).
\end{equation*}

By using the conditions for a maximization of this Hamiltonian, we
get after some algebras: $I_t^R$, $x_{t+1}$ and $P_{t+1}^R$.

Now, the manufacturer is faced with the optimal control problem.
To obtain the Stackelberg strategy of the manufacturer, we
maximize the objective function of the manufacturer by its
Hamiltonian function. We fix the value of $I_t^R$, then get the
Hamiltonian function of the manufacturer

\begin{equation*}\label{a4}
\begin{split}
H_{t}^M&=J_t^M+P_{t+1}^M(x_{t+1})+u_t(P_t^R).
\end{split}
\end{equation*}

So, we obtain, $I_{t}^M$, $x_{t+1}$, $P_{t}^M$ and $u_{t+1}$.

In this level, the game is between supplier as the leader and the
manufacturer as the follower. In fact, the reaction functions of
two followers (retailer and manufacturer) are placed on the
objective function of the leader (supplier). We can find final
equilibrium point.

For fixed $I_t^M$ Hamiltonian function of the supplier is defined
by
\begin{equation*}
\begin{split}
H_t^{S}&=J_t^{S}+P_{t+1}^{S}(x_{t+1})+u_t'(P_t^M).
\end{split}
\end{equation*}
We place reaction functions of followers into the leader's
Hamiltonian function and we obtain

$I_t^S$, $P_{t+1}^{S}$, $u'_{t}(P_t^M)$ and $x_{t+1}$.

For solving the above optimal control problem, we chose an
algorithms given by Medanic and Radojevic which is an augmented
discrete Hamiltonian matrix \cite{medanic}.

First, we assume
\begin{equation*}
    \left[ {\begin{array}{c}
   \widetilde{x}_{t+1}\\
   \widetilde{P}_t\\
  \end{array} } \right]=\left[ {\begin{array}{cccc}
   A&B\\
   C&D\\
  \end{array} } \right]\left[ {\begin{array}{c}
   \widetilde{x}_{t}\\
   \widetilde{P}_{t+1}\\
  \end{array} } \right]+\left[ {\begin{array}{c}
   D\\
   E\\
  \end{array} } \right]=\left[ {\begin{array}{c}
   A \widetilde{x}_t+B \widetilde{P}_{t+1}+D\\
   C \widetilde{x}_t+A \widetilde{P}_{t+1}+E\\
  \end{array} } \right],
\end{equation*}
where $
  \widetilde{x}_{t+1}=
  \left[ {\begin{array}{c}
   x_{t+1}\\
   u_{t+1}\\
  \end{array} } \right]
$ \quad and\quad $
  \widetilde{P}_{t+1}=
  \left[ {\begin{array}{c}
   p_{t+1}^M\\
   p_{t+1}^R\\
  \end{array} } \right]$,

$A, B,$ and $C$ are $2\times2$ matrices, and $D$ and $C$ are
$2\times1$ matrices.

We solved the above problem by sweep method \cite{Bryson}, by
assuming a linear relation between $\widetilde{p_t}$ and
$\widetilde{x_t}$; thus, we can compute the value of
$\widetilde{p_t}$ and $\widetilde{x}_{t}$ and we can obtain the
values of the other variables for all points in time by backward
and forward loop.

\section{Conclusion}

This paper investigated a decentralized three-tier supply chain
consisting of supplier, manufacturer and retailer for the
allocating CSR to members of supply chain system in over time. We
considered a two-level Stackelberg game consisting of two
followers and one leader. The members of a supply chain play games
with each other to maximize their own profits; thus, the model
used to be a long-term co-investment game model. The equilibrium
point at which members make their decisions to maximize profits by
implementing CSR among members of the supply chain in a time
horizon was determined. We applied control theory and used an
algorithm (augmented discrete Hamiltonian matrix) to obtain an
optimal solution for the dynamic game model.

\end{document}